\documentclass[11pt, twoside]{article}
\usepackage{amssymb, amsmath}
\usepackage[matrix,arrow]{xy}
\usepackage{amscd}

\begin{document}
\newtheorem{Mthm}{Main Theorem.}
\newtheorem{Thm}{Theorem}[section]
\newtheorem{Prop}[Thm]{Proposition}
\newtheorem{Lem}[Thm]{Lemma}
\newtheorem{Cor}[Thm]{Corollary}
\newtheorem{Def}[Thm]{Definition}
\newtheorem{Guess}[Thm]{Conjecture}
\newtheorem{Ex}[Thm]{Example}
\newtheorem{Rmk}[Thm]{Remark}
\newtheorem{Not}{Notation.}

\renewcommand{\theThm} {\thesection.\arabic{Thm}}
\renewcommand{\theProp}{\thesection.\arabic{Prop}}
\renewcommand{\theLem}{\thesection.\arabic{Lem}}
\renewcommand{\theCor}{\thesection.\arabic{Cor}}
\renewcommand{\theDef}{\thesection.\arabic{Def}}
\renewcommand{\theGuess}{\thesection.\arabic{Guess}}
\renewcommand{\theEx}{\thesection.\arabic{Ex}}
\renewcommand{\theRmk}{\thesection.\arabic{Rmk}}
\renewcommand{\theMthm}{}
\renewcommand{\theNot}{}
\newcommand{\spec}{\textnormal{Spec}\hspace{1mm}}
\newcommand{\proj}{\textnormal{Proj}\hspace{1mm}}
\newcommand{\Diff}{\textnormal{Diff}\hspace{1mm}}
\newcommand{\pf}{{\bfseries\itshape Proof. }}
\newcommand{\mult}{\textnormal{mult}}
\newcommand{\com}{\hspace{-2mm}\textnormal{\textbf{.}}\hspace{2mm}}
\newcommand{\bir}{-\hspace{-1mm}\rightarrow}
\newcommand{\qed}{\hfill Q.E.D.\newline}

\title{Birational Maps of Del Pezzo Fibrations.
}
\author{JIHUN PARK}
\date{}
\maketitle
\thanks{\centering Department of Mathematics\\
The Johns Hopkins University\\
Baltimore, Maryland, USA.\\}
\thispagestyle{empty}
\renewcommand{\thesection}{\large{\arabic{section}.}}
\section{\hspace{-3mm}\large{Introduction.}}
\renewcommand{\thesection}{\arabic{section}}
In classical result, it is known that 
any $\mathbb{P}^1$-bundle over a nonsingular complex curve
$T$ can be birationally transformed to a 
$\mathbb{P}^1$-bundle over 
$T$ by an elementary transformation.
Here, we can ask if it is also possible in 3-fold case. In other words, 
is it true that any nonsingular del Pezzo fibration over a nonsingular curve
can be transformed to another nonsingular del Pezzo fibration? 
In this question, we can add more condition on del Pezzo fibrations with
some kind of analogue from ruled surface cases, that is,
we can assume that their fibers are always nonsingular even though this 
is not true for any nonsingular del Pezzo fibration.

We ask the same question for local cases. Of course, we can 
birationally transform any $\mathbb{P}^1$-bundle over a germ of nonsingular 
complex curve $(T,o)$ into another $\mathbb{P}^1$-bundle over $(T,o)$.
But, in del Pezzo fibrations over $(T,o)$, something different
happens. In this paper, we will show that any del Pezzo fibration
of degree $d\leq 4$ with nonsingular special fiber cannot be birationally
transformed into another del Pezzo fibration with nonsingular special fiber.

Let $\mathcal{O}$ be a discrete valuation ring such that its residue field
$k$ is of characteristic zero. We denote $K$ the quotient field of 
$\mathcal{O}$.
Let $X_K$ be a variety defined over $\spec K$.
A model of $X_K$ is a flat scheme $X$
defined over $\spec \mathcal{O}$ whose generic fiber is isomorphic to $X_K$.
Fano fibrations are models of
nonsingular Fano variety defined over $K$.
In particular, del Pezzo fibrations of degree $d$ are models of 
nonsingular del Pezzo surfaces of degree $d$ defined over
$K$.
Del Pezzo fibrations are studied in \cite{Co} and \cite{Ko2}.
They constructed ``standard model" (\cite{Co})
and ``semistable model" (\cite{Ko2}) in each paper.

Now, we state the theorem which we will prove in this paper.

\begin{Mthm}
Let $X$ and $Y$ be del Pezzo fibrations of degree $d\leq 4$ over
$\spec \mathcal{O}$. Suppose that each scheme-theoretic special fiber is
nonsingular. Then any birational map between $X$ and $Y$  over
$\spec \mathcal{O}$ which is identical over generic fiber
is a biregular morphism.
\end{Mthm}

We should remark here that
even though it is hard to find such
examples, there are
del Pezzo fibrations
of degree $d\leq 4$ over $\spec \mathcal{O}$ with nonsingular special fibers
which can be birationally transformed
into another del Pezzo fibration
over $\spec \mathcal{O}$ with reduced and irreducible special fiber.
But, as in Minimal model program over 3-folds, we have to allow some
mild singularities, such as terminal ones, on them. In the end of this
paper, we will give such examples.

From now on, we explain standard definitions and notations
for this paper.
For more detail,
we can refer to  \cite{Ko1}, \cite{Pro99},
and \cite{Sh93}.

A variety $X$ means an integral scheme of finite type over a fixed field
$k$.
A log pair $(X,B)$ is a normal variety $X$ equipped with
a $\mathbb{Q}$-Weil divisor $B$ such that $K_X+B$ is $\mathbb{Q}$-Cartier.
A log variety is a log pair $(X,B)$ such that $B$ is a subboundary.

The discrepancy of a divisor $E$ over $X$ with respect to a log pair
$(X, B)$ will be denoted by $a(E;X,B)$.
And we will use the standard abbreviation plt, klt, and lc for
purely log terminal, Kawamata log terminal, and log canonical, respectively.

Let $(X,B)$ be an lc pair and $D$ an effective $\mathbb{Q}$-Cartier divisor
on $X$. The log canonical threshold (or lc threshold) of $D$ is the number
\[lct(X,B,D):=\sup\{c\phantom{s}|\phantom{s}
(X,B+cD) \phantom{s}is\phantom{s} lc\}.\]
If $B=0$, then we use $lct(X,D)$ instead of $lct(X,0,D)$.

Finally, we will use V. V. Shokurov's 1-complement which is a main tool
for this paper.
Let $X$ be a normal variety
and let $D$ be a reduced and irreducible divisor 
on $X$.
A divisor $K_X+D$, not necessarily log canonical, is $1$-complementary
if there is an integral Weil divisor $D^+$ such that
$K_X+D^+$ is linearly trivial,
$K_X+D^+$ is lc, and
$D^+\geq D$.
The divisor $K_X+D^+$ is called a 1-complement of $K_X+D$.
This is just special case of $n$-complements. But, it is enough for this paper.
For more detail about complements, we can refer to \cite{Pro99}, \cite{Sh93},
or \cite{Sh97}.
\vspace{1cm}\newline
\textbf{Acknowledgments.}  The author would like to thank Prof. V. V. Shokurov
for his invaluable support.

\renewcommand{\thesection}{\large{\arabic{section}.}}
\section{\hspace{-3mm}\large{Properties of certain birational maps.}}
\renewcommand{\thesection}{\arabic{section}}
Let $\mathcal{O}$ be  a discrete valuation ring
with local parameter $t$. The quotient field and residue field of
$\mathcal{O}$ are denoted by $K$ and $k$, respectively.
We always assume that the field $k$ is
of characteristic zero.
We denote $T=\textnormal{Spec}\hspace{1mm}\mathcal{O}$.
For a scheme $\pi:Z\longrightarrow T$,
its scheme-theoretic special fiber $\pi^*(o)$ is denoted by $S_Z$, where
$o$ is the closed point of $T$.
From now on, a birational map is always assumed to be identical
when restricted to the generic fibers.

Let $X/T$  be a $\mathbb{Q}$-factorial 
Gorenstein model of a nonsingular variety defined over $K$ which satisfies
the following three conditions.
\begin{itemize}
\item
(Special fiber condition)\\
The special fiber $S_X$ is a reduced and irreducible variety with nonempty
anticanonical linear system. 
Moreover,
log pair $(X,S_X)$ is plt.
\item
(1-complement condition)\\
For any $C\in|-K_{S_X}|$, there exists 1-complement
$K_{S_X}+C_X$ of $K_{S_X}$ such that $C_X$ does not contain any center
of log canonicity of $K_{S_X}+C$.
\item
(Surjectivity condition)\\
Any 1-complement of $K_{S_X}$  
can be extended to a 1-complement of
$K_X+S_X$. 
\end{itemize}
With Special fiber condition, we can easily show that $X$ has at worst
terminal singularities. Moreover, the special fiber $S_X$  is a 
variety over $k$ with Gorenstein canonical singularities.

Let $\phi:X\bir Y$ be a birational map over $T$, where
$X$ and $Y$ are $\mathbb{Q}$-factorial
Gorenstein models of a nonsingular variety defined over $K$ which satisfy
above three conditions.
Suppose that $\phi:X\bir Y$ is not an isomorphism in codimension 1. 
We fix a resolution of indeterminacy of $\phi:X\bir Y$ as follows.
\[ \xymatrix{
    & W\ar[dl]_f\ar[dr]^g\\
   X \ar@{-->}[rr]^{\phi}  \ar[dr] && Y \ar[dl]\\
  &T
}\]
Let $\widetilde{S_X}$ and $\widetilde{S_Y}$ be proper transformations
of $S_X$ and $S_Y$ by $f$ and $g$, respectively.
Since birational map $\phi$ is not an isomorphism in codimension 1,
$\widetilde{S_X}$ is a $g$-exceptional divisor and $\widetilde{S_Y}$
is  $f$-exceptional.

\begin{Lem}\label{shokurov}\com
Let $K_X+S_X+D_X$ be a 1-complement of $K_X+S_X$. And let $D_Y=\phi_*D_X$.
For any prime divisor $E$ over $X$,
\[a(E;X,qS_X+D_X)=a(E;Y,\alpha_q S_Y+D_Y),\]
where $q$ is any given number and $\alpha_q=-a(\widetilde{S_Y};X,qS_X+D_X)$.
Moreover, log canonical divisor $K_Y+S_Y+D_Y$ is linearly trivial.
\end{Lem}
\pf
Suppose that $E$ is a divisor on $W$.
Note that $f^{-1}_*D_X=g^{-1}_*D_Y=D_W$.
Then we have
\[K_W+q\widetilde{S_X}+D_W=
f^*(K_X+qS_X+D_X)-\alpha_q \widetilde{S_Y}+\sum a_iE_i,\]
and
\[K_W+\alpha_q \widetilde{S_Y} + D_W=g^*(K_Y+\alpha_q S_Y +D_Y)
+b\widetilde{S_X}+\sum b_iE_i,\]
where each $E_i$ is $f$-exceptional and $g$-exceptional.
From them, we get
\[f^*(K_X+qS_X+D_X)-g^*(K_Y+\alpha_q S_Y +D_Y)
=(q+b)\widetilde{S_X}+\sum (b_i-a_i)E_i.\]
Since $K_X+qS_X+D_X$ is numerically trivial,
we have
\[(q+b)\widetilde{S_X}+\sum (b_i-a_i)E_i\equiv_g 0.\]
By Negativity lemma, $b=-q$ and $b_i=a_i$. This prove the first statement.

Since $\phi$ is identical on generic fiber, it is clear that $D_Y$ is
linearly equivalent to $-K_Y$. 
Thus, the second statement follows from the fact that
$S_Y$ is linearly trivial. 
\qed
\begin{Lem}\label{wlog1}\com
There exists 1-complement $K_{S_X}+C_X$ (resp. $K_{S_Y}+C_Y$) of
$K_{S_X}$ (resp. $K_{S_Y}$) does not contain the center of 
$S_Y$ (resp. $S_X$) on $X$ (resp. $Y$).
\end{Lem}
\pf
Let $K_Y+S_Y+L_Y$ be a 1-complement of $K_Y+S_Y$.
By lemma~\ref{shokurov}, $a(\widetilde{S_Y};X,S_X+L_Y)\leq -1$,
where $L_X=\phi^{-1}_*L_Y$.
Clearly, the center of $\widetilde{S_Y}$ on $X$ is contained in
$C=L_X|_{S_X}$.
By inversion of adjunction, the center of $\widetilde{S_Y}$ is
a center of log canonicity singularities of $K_{S_X}+C$.
Furthermore, $K_{S_X}+C$ is linearly trivial by lemma~\ref{shokurov}.
Therefore, 1-complement condition implies the statement. 
\qed
\begin{Lem}\label{wlog2}\com
There is a 1-complement $K_X+S_X+D_X$ (resp. $K_Y+S_Y+H_Y$) of
$K_X+S_X$ (resp. $K_Y+S_Y$) such that  $D_X$ (resp. $H_Y$) 
does not contain the center
of $S_Y$ (resp. $S_X$).
\end{Lem}
\pf
It immediately follows from lemma~\ref{wlog1} and Surjectivity condition.
\qed

From now on, we fix 1-complements $K_X+S_X+D_X$ and $K_Y+S_Y+H_Y$
of $K_X+S_X$ and $K_Y+S_Y$, respectively, which satisfy the condition
in lemma~\ref{wlog2}.
We will use the notation $D_Y$, $D_W$, $H_X$ and $H_W$
for $\phi_*D_X$, $f_{*}^{-1}D_X$, $\phi_*^{-1}H_Y$
and $g_{*}^{-1}H_Y$, respectively.
Note that $g_{*}^{-1}D_Y=f_{*}^{-1}D_X$ and $g_{*}^{-1}H_Y=f_{*}^{-1}H_X$.

Now, we define the following condition.
\begin{itemize}
\item
(Total lc threshold condition)\\
The inequality $\tau_X+\tau_Y>1$ holds,
where 
$\tau_X=\min\{lct(S_X, C):C\in |-K_{S_X}|\}$ 
and
$\tau_Y=\min\{lct(S_Y, C):C\in |-K_{S_Y}|\}$.
\end{itemize}
\begin{Thm}\label{imissyou}\com
Under Total lc threshold condition, birational map $\phi$ is an isomorphism
in codimension 1.
\end{Thm}
\pf Suppose that $\phi$ is not an isomorphism in codimension 1.
We pay attention to the following eight equations;
\[K_W=f^*(K_X)+a\widetilde{S_Y}+\sum a_iE_i,\hspace{10mm}
\widetilde{S_X}=f^*(S_X)-b\widetilde{S_Y}-\sum b_iE_i,\]
\[D_W=f^*(D_X)-\sum c_i E_i,\hspace{10mm}
H_W=f^*(H_X)-e\widetilde{S_Y}-\sum e_iE_i,\]
\[K_W=g^*(K_Y)+n\widetilde{S_X}+\sum n_iE_i,\hspace{10mm}
\widetilde{S_Y}=g^*(S_Y)-m\widetilde{S_X}-\sum m_iE_i,\]
\[D_W=g^*(D_Y)-l\widetilde{S_X}-\sum l_iE_i,\hspace{10mm}
H_W=g^*(H_Y)-\sum r_iE_i.\]
First of all, $b=m=1$ since $S_X$ and $S_Y$ are reduced and irreducible.
Since $D_X$ does not contain the center of $\widetilde{S_Y}$ on $X$,
we have $\mult_{\widetilde{S_Y}}D_X=0$. For the same reason, we also have 
$\mult_{\widetilde{S_X}}H_Y=0$.  

By lemma~\ref{shokurov},
we get
$n+a-l=a(\widetilde{S_X};Y,-aS_Y+D_Y)=a(\widetilde{S_X};X,D_X)=0$ and
$a+n-e=a(\widetilde{S_Y};X,-nS_X+H_X)=a(\widetilde{S_Y};Y,H_Y)=0$, and hence
$a+n=l=e$. Since $X$ and $Y$ have at worst terminal singularities,
$a+n=l>0$.

Since $K_Y+S_Y+D_Y$ is linearly trivial by lemma~\ref{shokurov}, 
$(K_Y+S_Y+D_Y)|_{S_Y}=K_{S_Y}+D_Y|_{S_Y}$
is linearly trivial. 
Thus, $D_Y|_{S_Y} \in |-K_{S_Y}|$. 
Consequently,
it follows from inversion of adjunction
that $K_X+S_X+\tau_XH_X$ is lc.
By the same reason, $K_Y+S_Y+\tau_YD_Y$ is lc.

Now, we have $a(\widetilde{S_Y};X,S_X+\tau_XH_X)=a-1-\tau_Xe\geq -1$
and $a(\widetilde{S_X};Y,S_Y+\tau_YD_Y)=n-1-\tau_Yl\geq -1$.
But, $l=a+n\geq \tau_Xe+\tau_Yl=(\tau_X+\tau_Y)l >l$ 
by Total lc threshold condition.
Since $l>0$, this is impossible.
\qed

\renewcommand{\thesection}{\large{\arabic{section}.}}
\section{\hspace{-3mm}\large{Lc thresholds on
nonsingular del Pezzo surfaces.}}
\renewcommand{\thesection}{\arabic{section}}
Nonsingular del Pezzo surfaces were quite fully studied long time ago.
Furthermore, we understand singular del Pezzo surfaces very well. For example,
\cite{BW},
\cite{Dema80}, \cite{HW},
and \cite{R94} give us rich information.   
In this section, we will study some classical result on anticanonical
linear systems on del Pezzo surfaces with a modern point of view. 
Strictly speaking, we 
investigate all possible singular effective anticanonical
divisors on nonsingular del Pezzo surfaces. 
From this investigation, we can get some information 
on lc thresholds on nonsingular del Pezzo surfaces. 
\begin{Lem}\label{reduced}\com
Let $S$ be a nonsingular del Pezzo surface of degree $d\leq 4$.
Then, $K_S+C$ is lc in codimension 1 for any $C\in |-K_S|$.
\end{Lem}
\pf
Let $C=\sum^n_{i=1}m_iC_i\in |-K_S|$, where $C_i$'s are distinct
integral curves on $S$
and each $m_i\geq 1$.

First, we claim that if $C$ is not irreducible, then each $C_i$ is isomorphic
to $\mathbb{P}^1$. Suppose that $C_i$ is not isomorphic to $\mathbb{P}^1$.
Then, the self-intersection number of $C_i$ is greater than $0$.
Because $-K_S$ is ample, $C$ is connected. So, we have
\[2p_a(C_i)-2=(C_i+K_S)\cdot C_i=(1-m_i)C_i^2-\sum_{i\ne j}m_jC_j\cdot C_i<0,\]
which is contradiction. Thus, each component is a nonsingular rational curve.

Since $d=C\cdot(-K_S)=\sum^n_{i=1}m_iC_i\cdot(-K_S)$ and $-K_S$ is ample,
we have $\sum^n_{i=1}m_i\leq d$.

If $d=1$, then $n=1$ and $m_1=1$.

If $d=2$, then we have three possibilities $C_1$, $C_1+C_2$, and $2C_1$.
But the last case is absurd because the Fano index of $S$ is one.

Suppose $d=3$.
Then possibilities are $C_1$, $C_1+C_2$, $C_1+C_2+C_3$, $C_1+2C_2$,
$2C_1$, and $3C_1$. With the Fano index one, we can get rid of
the last two cases. For the case of $C=C_1+2C_2$, we consider the equation
$3=K_S^2=(C_1+2C_2)^2=C_1^2+4C_1\cdot C_2+4C_2^2$.
Since $(C_1+2C_2)\cdot(-K_S)=3$, we have $C_1\cdot(-K_S)=C_2\cdot(-K_S)=1$,
and hence $C_1^2=C_2^2=-1$. Thus, $C_1\cdot C_2=2$.
But, this implies contradiction $-2=2p_a(C_1)-2=C_1\cdot(-2C_2)=-4$.

Finally, we suppose that $d=4$.
We have eleven candidates,
$C_1$, $C_1+C_2$, $C_1+C_2+C_3$, $C_1+C_2+C_3+C_4$,
$C_1+2C_2$, $C_1+3C_2$, $C_1+C_2+2C_3$,
$2C_1+2C_2$, $2C_1$, $3C_1$, and $4C_1$.
Again, we can exclude the last four candidates by Fano index.
For the case of $C=C_1+3C_2$, we consider the equation
$4=K_S^2=(C_1+3C_2)^2=C_1^2+6C_1\cdot C_2+9C_2^2$. As before, we can see
$C_1^2=C_2^2=-1$. So, we have contradiction $3C_1\cdot C_2=7$.
Let's consider the case of $C=C_1+2C_2$. Since $(C_1+2C_2)\cdot(-K_S)=4$,
$C_1^2=0$ and $C_2^2=-1$. Then, we have $4=(C_1+2C_2)^2=-4+4C_1\cdot C_2$.
But, $-2=p_a(C_1)-2=-2C_1\cdot C_2$.
Finally, we consider $C=C_1+C_2+2C_3$. Then, each $C_i$ is $-1$-curve.
Since $4=(C_1+C_2+2C_3)^2=C_1^2+C_2^2+4C_3^2+2(C_1\cdot C_2+2C_1\cdot C_3
+2C_2\cdot C_3)$, we have $5=C_1\cdot C_2+2C_1\cdot C_3+2C_2\cdot C_3$.
But, $-2=2p_a(C_1)-2=-(C_2+2C_3)\cdot C_1$, and hence
$3=2C_2\cdot C_3$. But this is impossible.
\qed

Let $S$ be a nonsingular del Pezzo surface with Fano index $r$.
Then, there is an ample integral divisor $H$, 
called fundamental class of $S$, such that
$-K_S=rH$.
A curve $C$ on $S$ is called a line (resp. conic and cubic) 
if $C\cdot H=1$ (resp. 2 and 3).
\begin{Prop}\label{oum}\com
Let $S$ be a nonsingular del Pezzo surface of degree $d\leq 4$ and let
$C\in |-K_S|$. Suppose that $K_S+C$ is worse than lc.
\begin{enumerate}
\item
If $d=1$, then $C$ is a cuspidal rational curve.
\item
If $d=2$, then $C$ is one of the following;
\begin{itemize}
\item
$C=C_1+C_2$, where $C_1$ and $C_2$ are lines intersecting tangentially
at one point with $C_1\cdot C_2=2$.
\item
$C$ is a cuspidal rational curve.
\end{itemize}
\item
If $d=3$, then $C$ is one of the following;
\begin{itemize}
\item
$C=C_1+C_2+C_3$, where $C_1$, $C_2$, and $C_3$ are lines
intersecting at one point with $C_1\cdot C_2=C_1\cdot C_3=C_2\cdot C_3=1$.
\item
$C=C_1+C_2$, where $C_1$ and $C_2$ are a line and a conic 
intersecting tangentially at one point with $C_1\cdot C_2=2$.
\item
$C$ is a cuspidal rational curve.
\end{itemize}
\item
If $d=4$, then $C$ is one of the following;
\begin{itemize}
\item
$C=C_1+C_2+C_3$, where $C_1$ and  $C_2$ are lines, and $C_3$ is a conic 
intersecting at one point with $C_1\cdot C_2=C_1\cdot C_3=C_2\cdot C_3=1$.
\item
$C=C_1+C_2$, where $C_1$ and $C_2$ are a line and a cubic 
intersecting tangentially at one point with $C_1\cdot C_2=2$.
\item
$C=C_1+C_2$, where $C_1$ and $C_2$ are conics 
intersecting tangentially at one point with $C_1\cdot C_2=2$.
\item
$C$ is a cuspidal rational curve.
\end{itemize}
\end{enumerate}
\end{Prop}
\pf 
Note that
if $C$ is irreducible, then arithmetic genus $p_a(C)$ of $C$ is one.
If $C$ is not irreducible, then each component is isomorphic to
$\mathbb{P}^1$. And we can see the intersection numbers of two
different components of $C$ are less than or equal to $2$.

We can easily check the cases of degree 1 and 2.

Now, we suppose that $d=3$.
And we suppose that $C=C_1+C_2+C_3$.
Since $3=(C_1+C_2+C_3)\cdot(-K_S)$, each $C_i$ is a line.
From $2=2-2p_a(C_1)=C_1\cdot(C_2+C_3)$ and $3=C_1^2+C_2^2+C_3^2+
2C_1\cdot(C_2+C_3)+2C_2\cdot C_3$, we get $C_2\cdot C_3=1$.
Similarly, we can get $C_1\cdot C_2=C_1\cdot C_3=1$.
Since $K_S+C$ is not lc, these three lines intersect each other at one point.

If $C$ has less than 4 components, then we can show our statement
with the same method as above.

The only remaining that we have to show is that $K_S+C$ is lc if $d=4$ and
$C=C_1+C_2+C_3+C_4$. Since each $C_i$ is a line, we get
\[4=C^2=-4+2(C_1\cdot C_2+C_1\cdot C_3+C_1\cdot C_4
+C_2\cdot C_3+C_2\cdot C_4+C_3\cdot C_4).\]
And, we have
$C_1\cdot(C_2+C_3+C_4)=2-2p_a(C_1)=2$, $C_2\cdot(C_1+C_3+C_4)=2$,
$C_3\cdot(C_1+C_2+C_4)=2$, and $C_4\cdot(C_1+C_2+C_3)=2$.
With these 5 equations and connectedness of $C$, we can see that $C$ is
a normal crossing divisor. Thus, $K_S+C$ is lc.
\qed
\begin{Cor}\label{sunjoo}\com
Let $S$ be a nonsingular del Pezzo surface of degree $d\leq 4$.
\begin{itemize}
\item
If $d=1$, then $K_S+\frac{5}{6}C$ is lc for any $C\in |-K_S|$.
\item
If $d=2$, then $K_S+\frac{3}{4}C$ is lc for any $C\in |-K_S|$.
\item
If $d=3$ or $4$, then $K_S+\frac{2}{3}C$ is lc for any $C\in |-K_S|$.
\end{itemize}
\end{Cor}
\pf
If $C$ is three nonsingular curves intersecting each other at single point
transversally, then  $lct(X,C)=\frac{2}{3}$.
If $C=C_1+C_2$ where $C_i$'s are nonsingular curves  intersecting tangentially
with $C_1\cdot C_2=2$, then we have $lct(X,C)=\frac{3}{4}$.
For the case of a cuspidal rational curve, $lct(X,C)=\frac{5}{6}$.
Thus, our statement immediately follows from 
proposition~\ref{oum}.
\qed
\begin{Rmk}\com
Let $S$ be a nonsingular del Pezzo surface of degree $d$.
Then, we have the maximum number $r$ such that $K_S+rC$ is lc for any
$C\in|-K_S|$. It is easy to show that
such $r$ is $\frac{1}{3}$ (resp. $\frac{1}{2}$)
if $d=9$, $7$, or $d=8$ and Fano index $1$ 
(resp. $d=5$, $6$ or $d=8$ and Fano index $2$).
\end{Rmk}

\begin{Rmk}\label{base}\com
If $S$ be a nonsingular del Pezzo surface of degree 1, then
$|-K_S|$ has exactly one base point.
We can easily check that any element in $|-K_S|$ is nonsingular
at this point.
\end{Rmk}

\renewcommand{\thesection}{\large{\arabic{section}.}}
\section{\hspace{-3mm}\large{Proof of main theorem.}}
\renewcommand{\thesection}{\arabic{section}}
In this section, we will use the same notations as in the second section.
\\
{\bfseries\itshape Proof of main theorem. }
Since $-K_X$ and $-K_Y$ are ample over $T$,
Surjectivity condition follows from \cite{Pro99a}. 
By the same reason, birational map $\phi$ cannot be an isomorphism in
codimension 1 unless it is biregular (see \cite{Co1}).

It is enough to check 1-complement condition and Total lc threshold condition
by theorem~\ref{imissyou}. 
Total lc threshold condition immediately follows from corollary~\ref{sunjoo}.
If $2\leq d\leq 4$, then it is clear that 1-complement condition holds.
In the case of degree 1, 1-complement condition can be derived from
remark~\ref{base}.
\qed
\begin{Cor}\com
Let $X$ be a del Pezzo fibration over $T$ of degree $\leq 4$ with nonsingular
scheme-theoretic special fiber.
Then, the birational automorphism group of $X/T$ is
the same as the biregular automorphism group of $X/T$.
\end{Cor}
\pf
Note that we always assume that birational map is identical on generic fiber.
The statement immediately follows from  the main theorem.
\qed
\\
As an easy application of theorem~\ref{imissyou}, 
we can get the following well-known
example.
\begin{Ex}\com
\textnormal{Let $Z$ be a $\mathbb{P}^1$-bundle over $T$. 
Suppose that the special fiber $S_Z$ has no $k$-rational point.
In particular, the residue field $k$ is not algebraically closed.
Then, there is no birational transform of $Z$ into another
$\mathbb{P}^1$-bundle over $T$, because the special fiber $S_Z$
satisfies Total lc condition. If $S_Z$ has a $k$-rational point,
then Total lc condition fails. Moreover, it can be birationally transformed
into another $\mathbb{P}^1$-bundle over $T$ by elementary transformations.
}
\end{Ex}
\renewcommand{\thesection}{\large{\arabic{section}.}}
\section{\hspace{-3mm}\large{Examples.}}
\renewcommand{\thesection}{\arabic{section}}
If we allow some mild singularities on del Pezzo fibrations, 
then we can find birational maps of del Pezzo fibrations over $T$ with
reduced and irreducible special fiber.
In each example, note
that one of two del Pezzo fibrations has terminal singularities.
Before taking examples, we will state easy
lemma which helps us to understand our examples.
\begin{Lem}\label{kuwata}\com
Let $f(x_1,\cdots,x_m,y_1,\cdots,y_n)=g(x_1,\cdots,x_m)+h(y_1,\cdots,y_n)$ be
a holomorphic function near $0\in \mathbb{C}^{m+n}$ and
let $D_f=(f=0)$ on $\mathbb{C}^{m+n}$, $D_g=(g=0)$ on $\mathbb{C}^{m}$, and
$D_h=(h=0)$ on $\mathbb{C}^{n}$.
Then
\[lct(\mathbb{C}^{m+n},D_f) =
min\{1, lct(\mathbb{C}^{m},D_g)+lct(\mathbb{C}^{n},D_h)\}.\]
\end{Lem}
\pf See \cite{Kuwata}.
\qed

\begin{Ex}\com
\textnormal{This example comes from \cite{Co} and \cite{Ko2}.
Let $X$ and $Y$ be subschemes of $\mathbb{P}^3_{\mathcal{O}}$
defined by equations
$x^3+y^3+z^2w+w^3=0$ and $x^3+y^3+z^2w+t^{6n}w^3=0$, respectively, where
$n$ is a positive integer.
Note that $X$ is nonsingular and $Y$ has single singular point of type $cD_4$
at $p=[0,0,0,1]$.
Then, we have a birational
map $\rho_n$ of $X$ into $Y$ defined by
$\rho_n([x,y,z,w])=[t^{2n}x,t^{2n}y,t^{3n}z,w]$.
Now, we consider a divisor $D\in|-K_X|$ defined by $z=w$.
This divisor $D$ has a sort of good divisor because $K_X+S_X+D$
is lc and $D\vert_{S_X}$ is a nonsingular elliptic curve on $S_X$.
But, the birational transform $\rho_{n*}(D)$ of $D$ by $\rho_n$ is worse
than before.
First, $\rho_{n*}(D)\vert_{S_Y}$ is three lines intersecting each other
at single point (Eckardt point) transversally on $S_Y$.
Furthermore,
we can see that $\rho_{n*}(D)$  on $Y$ is defined by $z=t^{3n}w$.
And, the log canonical threshold of $\rho_{n*}(D)$ is $\frac{4n+1}{6n}$ by
lemma~\ref{kuwata},
and hence $K_Y+\rho_{n*}(D)$ cannot be lc.
}
\end{Ex}
\begin{Ex}\label{cusp}\com
\textnormal{Let $Z$ and $W$ be subschemes of $\mathbb{P}^3_{\mathcal{O}}$
defined by equations
$x^3+y^2z+z^2w+t^{12m}w^3=0$ and $x^3+y^2z+z^2w+w^3=0$, respectively, where
$m$ is a positive integer.
Here, $Z$ has a singular point of type $cE_6$ 
at $[0,0,0,1]$ and $W$ is nonsingular.
We have a birational map $\psi_m$ of $Z$ into $W$ defined by
$\psi_m([x,y,z,w])=[t^{2m}x,t^{3m}y,z,t^{6m}w]$.
Again, we consider a divisor $H\in|-K_Z|$ defined by $z=w$.
For the same reason as above, $H$ is a good divisor.
But, the log canonical threshold of the birational transform 
$\psi_{m*}(H)$ of $H$ by $\psi_m$ is $\frac{5m+1}{6m}$.
Therefore, if $m>1$, then $K_W+\psi_{m*}(H)$ cannot be lc.
Note that $\psi_{m*}(H)\vert_{S_W}$ is a cuspidal rational curve on $S_W$.
}
\end{Ex}
\begin{Ex}\com
\textnormal{We consider birational map $\varphi_m=\psi_m^{-1}$ from $W$ to $Z$,
where $W$, $Z$ and $\psi_m$ are the same as in example~\ref{cusp}.
And, we pay attention to nonsingular
divisor $L\in|-K_W|$ on $W$ defined by $x=0$.
Then, we can see that $\varphi_{m*}(L)\vert_{S_Z}$ consists of
a line and a conic intersecting tangentially each other.
And, the log canonical threshold of $\varphi_{m*}(L)$ is $\frac{9m+1}{12m}$,
and hence $K_Z+\varphi_{m*}(L)$ is not lc.
}
\end{Ex}
The following two examples were constructed by M.~Grinenko. One is a del Pezzo
fibration of degree 2, and the other is of degree 1.
\begin{Ex}\com
\textnormal{Let $X$ and $Y$ be subschemes of 
$\mathbb{P}^3_{\mathcal{O}}(1,1,1,2)$
defined by equations $w^2+x^3y+x^2yz+z^4+t^4xy^3=0$ 
and $w^2+x^3y+xy^3+x^2yz+t^2z^4=0$,
respectively, where $w$ is of weight 2.
The map $\phi:X\bir Y$ defined by
$\phi(x,y,z,w)=(x,t^2y,z,tw)$ is birational.
Subscheme $X$ has a singular point of type $cD_5$ at $[0,1,0,0]$.
Subscheme $Y$ has two singular points of types $cD_6$ and $cA_1$
at $[0,0,1,0]$ and $[1,0,-1,0]$, respectively.}
\end{Ex}
\begin{Ex}\com
\textnormal{Let $Z$ be a subscheme of $\mathbb{P}^3_{\mathcal{O}}(1,1,2,3)$
defined by equation $w^2+z^3+xy^5+t^4x^5y=0$, where $z$ and $w$ are of weight
2 and 3, respectively.
Then, we have a birational automorphism $\alpha$ of $Z$ defined by
$\alpha(x,y,z,w)=(y,t^2x,t^2z,t^3w)$. Note that 
$Z$ has a singular point of type $cE_8$ at $[1,0,0,0]$.
}
\end{Ex}
\footnotesize
\nocite{Pro99}
\nocite{Ko1}
\nocite{Ko2}
\nocite{Sh93}
\nocite{Co1}
\nocite{Co}
\nocite{HW}
\nocite{Sh97}
\nocite{Pro99a}
\nocite{BW}
\nocite{Dema80}
\nocite{R94}
\nocite{Kuwata}
\bibliographystyle{alpha}
\bibliography{bir}
\vspace{5mm}
E-mail address : \hspace{1mm}jhpark@chow.mat.jhu.edu

\end{document}